\documentclass[11pt]{article}
\usepackage[a4paper, total={174mm,241.3mm}]{geometry}

\usepackage[utf8]{inputenc}
\usepackage{parskip}
\usepackage{float}
\usepackage[title]{appendix}
\usepackage{authblk}
\usepackage{hyperref}

\usepackage{amsmath}
\usepackage{amssymb}
\usepackage{bbm}
\numberwithin{equation}{section}
\usepackage{cases}

\usepackage{array}
\newcolumntype{\$}{>{\global\let\currentrowstyle\relax}}
\newcolumntype{^}{>{\currentrowstyle}}

\usepackage{booktabs}
\usepackage{multirow}
\usepackage{caption}
\usepackage[font=footnotesize]{caption}
\usepackage{lipsum}

\usepackage{subcaption}
\usepackage{cleveref}
\usepackage[export]{adjustbox}

\usepackage{tikz}
\usetikzlibrary{patterns.meta, arrows.meta}
\usepackage{xcolor}
\usetikzlibrary{calc}
\usepackage{adjustbox}
\usepackage{lmodern}
\usetikzlibrary{intersections}
\tikzstyle{every picture}+=[font=\sffamily]

\definecolor{lime}{HTML}{A6CE39}
\DeclareRobustCommand{\orcidicon}{%
	\begin{tikzpicture}
	\draw[lime, fill=lime] (0,0) 
	circle [radius=0.16] 
	node[white] {{\fontfamily{qag}\selectfont \tiny ID}};
	\draw[white, fill=white] (-0.0625,0.095) 
	circle [radius=0.007];
	\end{tikzpicture}
	\hspace{-2mm}
}

\foreach \x in {A, ..., Z}{%
	\expandafter\xdef\csname orcid\x\endcsname{\noexpand\href{https://orcid.org/\csname orcidauthor\x\endcsname}{\noexpand\orcidicon}}
}

\usepackage[square,numbers]{natbib}
\usepackage{bibentry}
\title{Accelerated windowing for the crew rostering problem with machine learning}
\author[1*,3]{Philippe Racette\orcidA{}}
\author[2]{Frédéric Quesnel\orcidB{}}
\author[3,4]{Andrea Lodi\orcidC{}}
\author[1]{François Soumis\orcidD{}}
\affil[1]{Department of Mathematics and Industrial Engineering and GERAD, Polytechnique Montréal, Montréal, Québec, Canada}
\affil[2]{Department of Analytics, Operations and Information Technology, Université du Québec à Montréal, Montréal, Québec H2X 3X2, Canada}
\affil[3]{Canada Excellence Research Chair in Data Science for Real-Time Decision-Making, Polytechnique Montréal, Quebec H3C 3A7, Canada}
\affil[4]{Cornell Tech and Technion - IIT, New York City, New York, United States}
\date{*Corresponding author. E-mail: philippe.racette@gerad.ca}

\newcommand{\midarrow}{\tikz \draw[-{Triangle[width=3mm]}] (0,0) -- +(0.1,0);}
\newcommand{\whitearrow}{\tikz \draw[-{Stealth[width=3mm]}] (0,0) -- +(0.1,0);}
\newcommand{\forwardarrow}{\tikz \draw[-{Stealth}{Stealth}] (0,0) -- +(0.1,0);}
\newcommand{\crossarrow}{\tikz \draw[-{Rays[width=3mm]}] (0,0) -- +(0.1,0);}

\usepackage[labelfont=bf]{caption}
\DeclareCaptionLabelSeparator{none}{ }
\begin{document}

\maketitle

\begin{abstract}
\noindent The crew rostering problem (CRP) for pilots is a complex crew scheduling task assigning pairings, or sequences of flights starting and ending at the same airport, to pilots to create a monthly schedule. In this paper, we propose an innovative solution method for the CRP that uses a windowing approach. First, using a combination of machine learning (ML) and combinatorial optimisation (CO), we quickly generate an initial solution. The solution is obtained with a sequential assignment procedure (\textit{seqAsg}) based on a neural network trained by an evolutionary algorithm. Then, this initial solution is reoptimized using a branch-and-price algorithm that relies on a windowing scheme to quickly obtain a CRP solution. This windowing method consists of decomposing the optimization horizon into several overlapping windows, and then optimizing each one sequentially. Although windowing has been successfully used in other airline applications, it had never been implemented for the CRP, due to its large number of horizontal constraints involving the whole planning horizon. We test our approach on two large real-world instances, and show that our method is over ten times faster than the state-of-the-art branch-and-price CRP solver GENCOL while providing solutions on average less than 1\% away from optimality. We show that our windowing approach greatly benefits from being initialized with good-quality ML-based solutions. This is because the initial solution provides reliable information on the following windows, allowing the solver to better optimize the current one. For this reason, this approach outperforms other naive heuristics, including stand-alone ML or windowing.
\end{abstract}



\section{Introduction}
The crew rostering problem (CRP) is a form of crew assignment concerned with creating sets of schedules for pilots and flight attendants called rosters. Taking a set of pairings (i.e., sequences of flights and layovers starting and ending at the same airport) as input,  the CRP forms schedules by assigning them to pilots over a predetermined time horizon. This is done while taking into consideration crew preferences regarding specific flights, days off and special requests made by crew members. A number of airline and collective agreement regulations must also be observed.

Developing efficient methods to solve crew scheduling tasks in general is of prime importance to airline carriers, as problems such as the CRP take a long time to solve. One class of methods that has gained ground over the last few years is the use of machine learning (ML) techniques in the context of combinatorial optimization (CO) and operations research (OR) problems.  In this paper, we build upon the previous work in Racette \textit{et al.} \cite{racette2024} and develop a method that quickly generates good-quality CRP solutions. They propose a sequential assignment procedure named \textit{seqAsg}, an ML-based approach that generates complete rosters to facilitate schedule planning. Although those rosters could be used to provide insight to planners, they were not of a good enough quality to be used
in practice. We expand this contribution by showing that these rosters can be leveraged to obtain different near-optimal solutions produced faster than in the existing literature, and better than what would be obtained by different heuristics.

To do so, we propose a windowing approach that can take ML-generated rosters as input and improve them with a branch-and-price algorithm. Windowing is a solution method that decomposes the horizon into several overlapping time windows, and sequentially optimizes each window. Although windowing was applied to great success to the crew pairing problem (CPP) such as in Saddoune, Desaulniers and Soumis \cite{saddoune2009}, it has not been attempted for the CRP. This is perhaps because the CPP has more localized horizontal constraints (each pairing lasting a few days). By contrast, the CRP has horizontal constraints spanning the whole horizon (schedule validity constraints). As the results of this paper show, windowing can cope with this difficulty well, and especially so if a good-enough initial solution is provided. This initial solution helps by providing an estimate of the horizontal contributions outside the current optimization window.

The benefits of windowing lie in the fact that by only enhancing part of a solution and freezing the rest, the size of the problem treated in each window decreases considerably. However, one possible drawback of windowing is the fact that the method is no longer optimal. The question is then to determine if any loss in the CRP objective is minimal enough to justify significant gains in computational time.

Given this context, we make the following contributions:
\begin{itemize}
   \item We propose a windowing approach to the CRP for pilots and prove its effectiveness in obtaining rosters within 1\% of optimality. Our method is over 10 times faster than the state of the art.
   \item We show that our windowing method greatly benefits from an initial solution, even if this solution is far from the optimum. Such a solution can be obtained using the ML-based method proposed by Racette \textit{et al.} \cite{racette2024}.
\end{itemize}
We also provide a brief discussion of the future implications of this work. In particular, we mention the application of fast ML methods to CO problems, and how this research adds support to the use of ML to enhance solving methods able to exploit even limited additional information.

We now give an outline of the paper. In Section \ref{revue}, we present a literature review of the topics relevant to this paper. In Section \ref{rostering}, we provide a detailed description of the CRP. In Section \ref{reopt}, we lay out the details of the windowing  procedure. In Section \ref{data_method}, we present our experimental protocol and provide some implementation details. In Section \ref{results}, we present the main results supporting our contributions. We conclude in Section \ref{discussion} with a discussion of the work presented in this paper and of some of its implications.

\section{Literature review} \label{revue}
In this section, we review the relevant literature for the work presented in this paper. We start with Section \ref{revcrp} by giving some context relevant to crew scheduling in general and the CRP more specifically. In Section \ref{tradition}, we provide more information about column generation, a technique often used to solve the CRP, and windowing. We describe the literature regarding solving OR problems with ML algorithms in Section \ref{ml_methods}. Finally, we provide the essential details related to the sequential assignment procedure we use in Section \ref{seqAsg}.

\subsection{Crew scheduling and the CRP} \label{revcrp}
Crew scheduling is a rich field that tackles many different problems. Gopalakrishnan and Johnson \cite{gopalakrishnan2005} describe the main tasks addressed by airline scheduling as a sequential process. The last two parts of the airline scheduling process are the CPP and the crew assignment problem, which includes crew bidding, crew rostering, and preferential bidding. As together these two steps make up crew scheduling, we describe them in more detail.

The CPP is concerned with building pairings, i.e., sequences of flights starting and ending at the same airport base, of minimal cost while serving all flights planned in the previous scheduling steps. Once pairings are made, they must be assigned to crew members, whether pilots or flight attendants, while also observing many additional, often complex regulations. Airlines have different ways to do so. One way is to have crew members bid on pre-made schedules, called bidlines, and to assign them while holding these bids in consideration. Another possibility is crew rostering, where an entire set of personalized schedules for the pilots is called a roster. In crew rostering, crew members state preferences over certain flights, days off and sometimes other special activities (e.g., training for flying certain kinds of aircrafts), and the schedule is optimized for a certain criterion while ensuring that every pairing is carried out in a way that respects all rules. This is the approach considered in this paper. A hybrid approach, called preferential bidding, creates schedules much as in the CRP, but by also taking seniority into account in the objective. We note that each of these crew scheduling problems may use different formulations with one key difference lying in the criteria taken into account by the model objective. A review of these formulations is found in Kohl and Karisch \cite{kohl2004}.

\subsection{Column generation and windowing} \label{tradition}
 Column generation is a technique that is often used to solve complex CO problems. It decomposes the resolution process in two alternating phases: the restricted master problem (RMP) phase, where a model considering only a subset of possible schedules is solved, and the subproblem phase, where new promising columns are generated and added to the next RMP. Column generation is often embedded into a branch-and-bound tree where branching steps occur throughout the solving process. This ensures that an integer solution is found. Such algorithms are called "branch-and-price" algorithms and are described in Barnhart \textit{et al.} \cite{barnhart1998} and Desrosiers \textit{et al.} \cite{desrosiers2024}. This method is now quite standard in commercial solvers.

While the branch-and-price framework presented above works well to solve many CO tasks in an efficient manner, there are some situations where we might want a faster solution than these methods can allow. Some examples include real-time optimization or planning which may have to be performed in a few seconds, or the reoptimization of a solution. One approach attempting to do this is windowing. With windowing, the planning horizon is broken into time windows of a few days that are treated separately. Windowing has been used in contexts as varied as flight and maintenance planning, the fishing industry, bus dispatching, the railway industry, and the CPP \cite{afsar2006, millar1998, gkiotsalitis2020, nielsen2011, saddoune2009}. In particular, Saddoune, Desaulniers and Soumis \cite{saddoune2009} have found that the resolution of the CPP is accelerated by an average of 23.1\% on 7 instances when compared to a sequential three-phase approach. However, windowing has not been attempted for the CRP, due to the presence of constraints extending for one month, or the entire horizon's length.

\subsection{ML, metaheuristics and crew scheduling} \label{ml_methods}
ML methods have gained traction in recent years and have been increasingly applied to help solving OR problems, as have metaheuristics. Both also often serve as comparison points for each other. For an in-depth description, see Bengio, Lodi and Prouvost \cite{bengio2021} and Scavuzzo \textit{et al.} \cite{scavuzzo2024}.

With respect to crew scheduling, several contributions have been made. Yaakoubi, Soumis and Lacoste-Julien \cite{yaakoubi2020a, yaakoubi2020b} use supervised learning to make clusters of flights that are likely to belong to the same pairing to speed up the resolution of the CPP with dynamic constraint aggregation. Pereira \textit{et al.} \cite{pereira2022} use supervised learning to outperform strong branching during the resolution of the CPP while using properties of the branching tree as features, a type of approach often called ``learning to branch". Quesnel \textit{et al.} \cite{quesnel2022} consider the CRP and learn which pairings are most likely to be part of a crew's schedule, reducing the size of the problem and improving computing time by a factor of up to 10. Similarly, sometimes only promising columns, arcs, or flight connections, respectively, are considered as part of the column generation algorithm to enhance computing time \cite{morabit2021, morabit2023, tahir2021}. Finally, Maenhout and Vanhoucke \cite{maenhout2010} solve the CRP with a scatter search heuristic for two-week time horizons, 100 pilots and 600 pairings and obtain solutions 1.67\% away from optimality in an average time of 127 seconds. Their work stands out in the literature attempting to solve the CRP with metaheuristics as the instances used are particularly large, and the solutions obtained are of good quality.

\subsection{Sequential assignment} \label{seqAsg}
In order to provide ML-generated solutions as input to the windowing implementation of our CRP solver, we use the sequential assignment procedure presented in Racette \textit{et al.} \cite{racette2024}. The \textit{seqAsg} procedure is a constructive heuristic that solves an assignment problem for each day of the CRP horizon (i.e., one month) in sequential order. Once this process is completed, a full roster, or set of schedules, is available. On each day, there is a set of pilots who are available to receive pairing assignments, single days off, or vacations. There is also a fixed set of pairings that must be assigned. Each assignment is given a utility (i.e., a measure of how desirable it is to make this specific assignment) intended to estimate the value of giving a pilot the activity specific to that assignment.

The utility function mentioned above must be learned. Racette \textit{et al.} \cite{racette2024} do this using reinforcement learning or RL. They learn a policy that gives an approximation of the value of an assignment given certain conditions in the learning environment. This environment includes a reward function that accounts for satisfaction and feasibility. The ML weights for the neural network used to approximate the utility of an assignment are learned with the CMA-ES evolutionary algorithm \cite{hansen2001}.

\section{Crew rostering problem} \label{rostering}
In this section, we describe the version of the CRP that we use in this paper. Since the data we present in Section \ref{Data} was originally found in Kasirzadeh, Saddoune and Soumis \cite{kasirzadeh2017}, the structure of the problem shown here follows essentially the same pattern. In particular, in this work, we consider the creation of personalized rosters, where pilots express preferences before scheduling starts and where these preferences are taken into account so as to maximize the overall crew satisfaction. Pilots declare weighted preferences, with a weight indicating the extent of the satisfaction given by its fulfilment. Each pilot receives a similar budget for the allocation of preference weights. Additionally, some pilots have prescheduled off period requests that must be honored.

Let $W$ be a set of pairings, $K$ a set of pilots, and $D$ the set of possible bases. Each pilot and each pairing is associated with one base $d \in D$. The pairings are spread over an horizon of exactly one month. The CRP is then concerned with assigning all pairings to the pilots while ensuring each pilot has a complete schedule and in such a way as to maximize the sum of satisfaction scores for the whole set of pilots.

Schedules are constrained by various rules and regulations. We have the following rostering restrictions. Each pairing must be assigned to exactly one pilot and each pilot included must receive a feasible schedule. Pilots may not fly more than $T^{flight}$ hours in a month, must not be assigned more than $T^{work}$ consecutive duties, i.e., calendar days with at least one flight departure, and must receive at least $T^{off}$ days off. All preassigned days off must be included, and, lastly, there must also be at least $T^{min}$ hours of uninterrupted rest between two consecutive pairings.

Concerning the objective, the main criterion is pilots' satisfaction with their schedules, which we define as follows. Pilot preferences fall into two categories: preferred flights (as opposed to pairings) and preferred three-day vacations. Then, let us consider pilot $k \in K$. Let $F_k$ be the set of preferred flights for pilot $k$ and let $f \in F_k$ be such a flight with weight $v_{f}^k$. Likewise, let $O_k$ be the set of preferred vacations for pilot $k$ and let $o \in O_k$ be such a vacation with weight $v_{o}^k$. Let $\Omega$ and $\Omega^k$ be the set of feasible schedules, and feasible schedules for pilot $k\in K$, respectively. Let $p \in \Omega^k$ be such a schedule. Let $a_{fp}^k$ be a constant with value 1 if flight $f$ is in one of the pairings in schedule $p$ assigned to pilot $k$, and 0 otherwise. Finally, let $a_{op}^k$ be a constant with value 1 if preferred vacation $o$ is in schedule $p$ assigned to pilot $k$, and 0 otherwise. Then, we can define the satisfaction score brought to pilot $k$ by schedule $p$, $c_p^{k}$, in the following way:
\begin{equation}
    c_p^{k} = \sum_{f \in F^k}a_{fp}^{k}v_{f}^{k} + \sum_{o \in O_k}a_{op}^{k}v_{o}^{k}. \label{eq:costs}
\end{equation}
As it can be seen from Equation \eqref{eq:costs}, a pilot's satisfaction with their schedule is the sum of all the weighted preferences granted to them. These scores are used in the objective formulation for the CRP model.

\subsection{Mathematical formulation} \label{math_eq}
In this section, we present the model we use to solve the CRP. Let $Q^k$ be the set of preassigned days off for pilot $k\in K$. Let $n_w, w \in W$ be the number of flights in pairing $w \in W$. Let $a_{wp}$ be a constant taking value 1 if pairing $w$ is in schedule $p \in \Omega^k, k \in K$, and 0 otherwise. Let $a_{qp}$ be a constant taking value 1 if preassigned day off $q \in Q^k$ is in schedule $p \in \Omega^k, k \in K$. Let $c_p^k$ be the cost (satisfaction score) of schedule $p \in \Omega^k, k \in K$, as defined by Equation \eqref{eq:costs}. Let $x_p^k$ be a variable equal to 1 if schedule $p\in \Omega^k$ is assigned to pilot $k$, and 0 otherwise. Finally, let $s_w$ and $y_q$ be the slack variables ensuring the coverage of pairing $w \in W$ and preassigned day off $q \in Q^k, k \in K$, respectively.

Also, within this model, it is possible for certain pairings not to be assigned at the cost of a penalty. Let this penalty be $C^{F}$ per unassigned flight in the pairing. In practice, such pairings would be assigned to reserve pilots, a situation judged undesirable due to the necessity of paying reserve pilots high wages for few hours worked. Additionally, to ensure feasibility, we allow not assigning prescheduled days off at the cost of a different penalty. Let its cost be $C^{D}$ for each such day. $C^{D}$ is set at a prohibitively high value so no preassigned days off are left unassigned unless absolutely necessary.

Then, the CRP reads as
\begin{equation}
    \max \sum_{k \in K}\sum_{p \in \Omega^{k}}c_{p}^{k}x_{p}^{k} - C^{F}\sum_{w \in W}n_{w}s_w - C^{D}\sum_{k \in K}\sum_{q \in Q^k}y_q \label{eq:crp:obj}
\end{equation}
subject to:
\begin{align}
\sum_{k \in K}\sum_{p \in \Omega^{k}}a_{wp}x_{p}^{k} + s_w = 1, & \quad \forall w \in W \label{eq:crp:cover} \\
\sum_{p \in \Omega^{k}}a_{qp}x_{p}^{k} + y_q = 1, & \quad  \forall k \in K, \forall q \in Q^{k} \label{eq:crp:off} \\
\sum_{p \in \Omega^{k}}x_{p}^{k} = 1, & \quad \forall k \in K \label{eq:crp:pilots} \\
x_{p}^{k} \in \{0,1\}, & \quad \forall k \in K, \forall p \in \Omega^{k} \label{eq:crp:bin}\\
s_w \in \{0,1\}, &  \quad \forall w \in W \label{eq:pt:bin} \\
y_q \in \{0,1\}, &  \quad \forall k \in K, \forall q \in Q^{k} \label{eq:pre:bin}
\end{align}

Objective function \eqref{eq:crp:obj} maximizes the sum of three terms: one that takes into account the total satisfaction provided to the whole set of pilots ready to work for that roster, and two more terms that add penalties to the objective for unassigned pairings and missing preassigned days off, respectively. Constraints \eqref{eq:crp:cover} ensure that each pairing is assigned. Constraints \eqref{eq:crp:off} ensure that each preassigned day off is assigned. Constraints \eqref{eq:crp:pilots} ensure that each pilot is given a schedule. Constraints \eqref{eq:crp:bin} ensure that assignment variables are binary (i.e., an assignment is selected or not). Constraints \eqref{eq:pt:bin} and \eqref{eq:pre:bin} ensure that all slack variables are binary. We note that constraints related to the number of consecutive work days, total flight time, number of days off and postpairing rest are not included as they can all be assessed on an individual schedule basis. By only generating schedules that meet these requirements at each iteration of column generation, we do not need to state these constraints explicitly in the model.

Typically, the number of elements in $\Omega$ is too large for them to be enumerated explicitly, and in practice, there is no set of feasible schedules readily available. In order to solve the problem efficiently, one approach is to embed the mathematical formulation stated in this section in a branch-and-price scheme, where linear relaxations of model \eqref{eq:crp:obj}-\eqref{eq:pre:bin} are solved with column generation inside a branch-and-bound tree.

\section{Methodology} \label{reopt}
In this section, we describe resolution by windowing in more depth. In Section \ref{subproblem}, we describe the structure of the subproblem networks that are used to generate new columns during column generation. In Section \ref{win_adapt}, we show how these networks are adapted to implement a windowing procedure. For a detailed review of the different kinds of optimization strategies used in crew scheduling and many of their examples, see Desaulniers, Desrosiers and Solomon \cite{desaulniers2002}.

\subsection{Subproblem networks} \label{subproblem}
To implement windowing as described in Section \ref{win_adapt}, it is essential to understand the subproblem phase of column generation mentioned in Section \ref{tradition}. As a reminder, column generation is divided in two alternating phases. During the RMP phase, the relaxation of problem \eqref{eq:crp:obj} - \eqref{eq:pre:bin} is solved with a restricted set of feasible schedules, named $\Omega'$. Next, the subproblem phase generates new schedules (i.e., columns) that are added to the restricted set of schedules considered when solving the RMP. For this to be possible, one needs a way of generating new feasible schedules and a method to check their reduced cost.

To compute the reduced cost of a candidate schedule, we use the dual solution of the RMP. Let us consider pilot $k \in K$. Let $\alpha_w$ denote the dual value for constraint \eqref{eq:crp:cover}  for pairing $w$. Let $\beta_{q}$ denote the dual value for constraint \eqref{eq:crp:off}  for pilot $k$'s preassigned day off $q$. Let $\gamma^{k}$ denote the dual value for constraint \eqref{eq:crp:pilots} for the scheduling requirement of pilot $k$.

The reduced cost of schedule $p$ for pilot $k \in K$ is then 
\begin{align} \label{eq:dual_gen}
    \overline{c}_{p}^{k} = {} & c_{p}^{k} - \sum_{w \in W}\alpha_{w}a_{wp}^{k} - \sum_{q \in Q^{k}}\beta_{q}a_{qp}^{k} - \gamma^{k}.
\end{align}

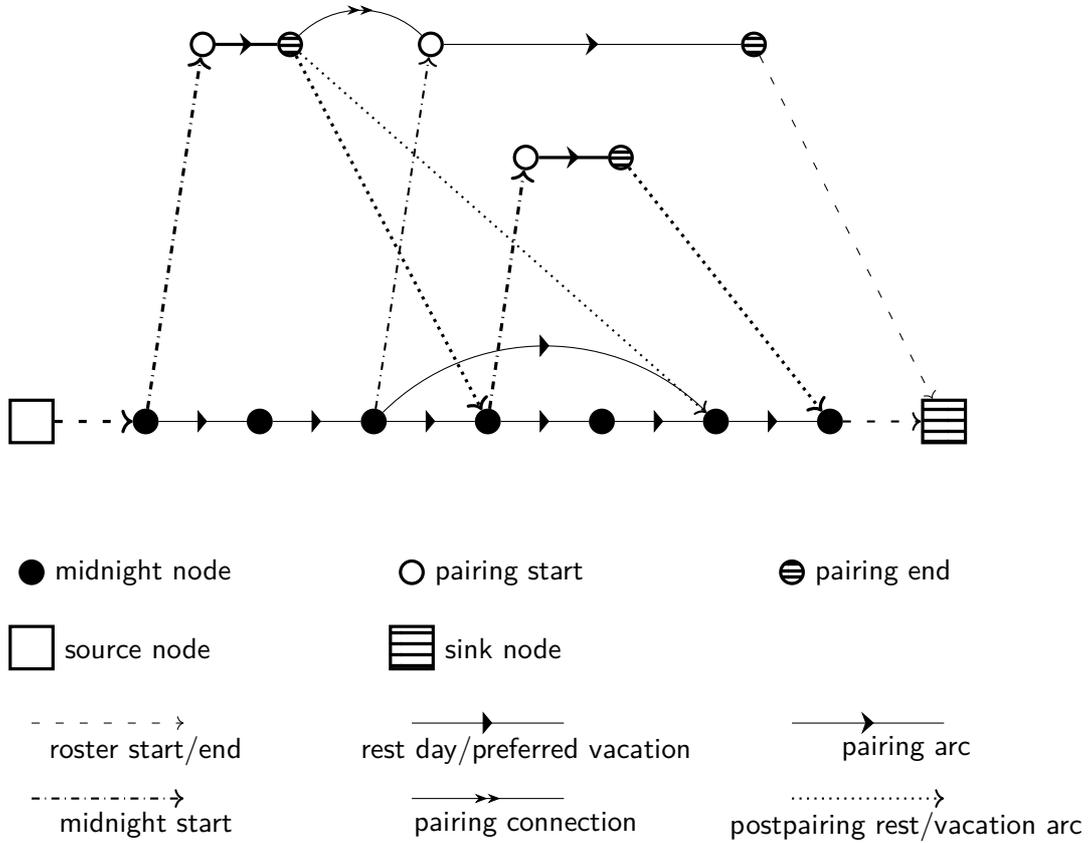
\begin{figure}[ht!]
\centering
\captionsetup{justification=centering, margin=1.75cm}
\begin{minipage}{\linewidth}
\centering
\begin{tikzpicture}
[start/.style={circle,draw=black, very thick, inner sep = 3pt},
finish/.style={circle,draw=black, very thick, inner sep = 3pt, pattern={Lines[angle=135,distance=2.5pt]}},
midnight/.style={circle,draw=black, very thick, fill=black, inner sep = 3pt},
source/.style={rectangle,draw=black, very thick, inner sep = 8pt},
sink/.style={rectangle, draw=black, very thick, inner sep = 8pt, pattern={Lines[angle=135,distance=4pt]}}]

\node at (-6,0) [source] (source0) {};
\node at (-3.75, 5) [start] (start01) {};
\node at (-4.5,0) [midnight] (midnight01) {};
\node at (-2.6, 5) [finish] (finish01) {};
\node at (-3,0) [midnight] (midnight02) {};
\node at (-1.5,0) [midnight] (midnight03) {};
\node at (0,0) [midnight] (midnight04) {};
\node at (-0.75, 5) [start] (start02) {};
\node at (1.5,0) [midnight] (midnight05) {};
\node at (3.5, 5) [finish] (finish02) {};
\node at (3,0) [midnight] (midnight06) {};
\node at (0.5,3.5) [start] (start03) {};
\node at (4.5,0) [midnight] (midnight07) {};
\node at (1.75,3.5) [finish] (finish03) {};
\node at (6,0) [sink] (sink0) {};

\draw [->, loosely dashed, very thick] (source0) -- (midnight01);
\draw (midnight01) -- node {\midarrow} (midnight02);
\draw (midnight02) -- node {\midarrow} (midnight03);
\draw (midnight03) -- node {\midarrow} (midnight04);
\draw (midnight04) -- node {\midarrow} (midnight05);
\draw (midnight05) -- node {\midarrow} (midnight06);
\draw (midnight06) -- node {\midarrow} (midnight07);
\draw [->, loosely dashed, thick] (midnight07) -- (sink0);

\draw[->, dash dot, very thick] (midnight01) -- (start01);
\draw[very thick] (start01) -- node {\whitearrow} (finish01);
\draw[->, dotted, very thick] (finish01) -- (midnight04);
\draw[->, dotted, thick] (finish01) -- (midnight06);

\draw[->, dash dot, thick] (midnight03) -- (start02);
\draw (start02) -- node {\whitearrow} (finish02);
\draw [->, loosely dashed] (finish02) -- (sink0);

\draw[->, dash dot, very thick] (midnight04) -- (start03);
\draw[very thick] (start03) -- node {\whitearrow} (finish03);
\draw[->, dotted, very thick] (finish03) -- (midnight07);

\draw (finish01) to[out=45, in=135] node {\forwardarrow} (start02);
\draw (midnight03) to[out=45, in=135] node {\midarrow} (midnight06);

\node at (-6, -2) [midnight, label=right:midnight node] {};
\node at (-1, -2) [start, label=right:pairing start] {};
\node at (4, -2) [finish, label=right:pairing end] {};
\node at (-6, -3) [source, label=right:source node] {};
\node at (-1, -3) [sink, label=right:sink node] {};

\draw [->, loosely dashed] (-6, -4) -- (-4, -4) node [near end, below=1pt] {roster start/end};
\draw (-1, -4) -- node {\midarrow} (1, -4) node [near end, below=1pt] {rest day/preferred vacation};
\draw (4, -4) -- node {\whitearrow} (6, -4) node [near end, below=1pt] {pairing arc};
\draw [->, dash dot, thick] (-6, -5) -- (-4, -5) node [near end, below=1pt] {midnight start};
\draw (-1,-5) -- node {\forwardarrow} (1,-5) node [near end, below=1pt] {pairing connection};
\draw[->, dotted, thick] (4,-5) -- (6,-5) node [near end, below=1pt] {postpairing rest/vacation arc};

\end{tikzpicture}
\end{minipage}

\caption{Example of subproblem network for a pilot}
\label{fig:roster}
\end{figure}

Since Equation \eqref{eq:dual_gen} is expressed as a sum of terms involving the total satisfaction provided to pilot $k \in K$ by schedule $p$, we model the subproblems as constrained shortest-path problems in an acyclic network. In each network, nodes represent different events such as the beginning or end of a pairing, and where arcs represent the selection of a certain action imposed to the pilot. Such a network is depicted in Figure~\ref{fig:roster}. The costs placed on the arcs of the network are then chosen so that the sum of the arc costs along a path from the source to the sink, representing a schedule, is equivalent to the reduced cost of the entire schedule. Paths are further constrained by resource windows on each node, which are described below.

We now describe the structure of the network for pilot $k$'s subproblem. The network includes five types of nodes. First, there is a \textit{source} node $o(k)$ for the pilot at the start of their schedule. There is also a \textit{sink} node $d(k)$ which marks the end of a schedule. Next, for each pairing in $W$, there is a \textit{beginning-of-pairing} node and an \textit{end-of-pairing} node. Finally, there is one \textit{midnight} node per day. These mark the beginning of a new calendar day. Every day off starts and ends with such a node.
 
Let $A^{k}$ be the arc set linking the nodes for in the subproblem of pilot $k \in K$. We note that we remove all arcs that make it possible to avoid assigning preassigned days off, although the associated penalties are still technically implemented in the used solver as a basic feature and remain part of objective \eqref{eq:crp:obj}. In that case, not granting pilot $k$'s preassigned days off is equivalent to not assigning them a schedule. This is done with a static column for that purpose.

$A^{k}$ includes the following types of arcs. First, there is a \textit{roster start} arc which represents the start of the schedule. It links the source node for pilot $k$ to the midnight node for the first day of the month. \textit{Roster end} arcs mark the end of the schedule. There is a roster end arc between the end node of each pairing ending on the last day of the horizon and the sink. There is also a roster end arc between the last midnight node of the month and the sink.
 
Two types of arcs represent days off or three-day vacations. \textit{Rest} arcs represent one single day off and link the midnight node marking the start of the day off to the midnight node marking the beginning of the next day. \textit{Vacation arcs}, on the other end, are used to assign three-day vacations. These arcs go from the midnight node at the beginning of the first day of the vacation to the midnight node three calendar days later.
 
\textit{Pairing} arcs represent carrying out pairings. The arc for a specific pairing links the pairing node at its beginning to the pairing node at its end. \textit{Midnight start} arcs are used to assign a pairing after one or several days off or on the first day of the horizon. They link the midnight node at the start of the day to the beginning-of-pairing node of the first pairing assigned on that day. \textit{Pairing connection} arcs ensure the direct transition between two pairings, one of which ends on a given day while the other starts on the same or following day. They go from the end-of-pairing node for the first pairing in the connection to the beginning-of-pairing node for the other pairing. Only those connection arcs that allow for sufficient postpairing rest time $T^{min}$ between pairings are included. 
 
Finally, some arcs are related to both pairing and rest periods. They ensure that time off is granted to the pilot following a pairing arrival. \textit{Postpairing rest} arcs are used to assigned a day off after a pairing arrival. The arc links the end-of-pairing node at the pairing arrival to the midnight node at the end of the next day. \textit{Postpairing vacation} arcs allow the assignment of a three-day vacation after a pairing arrival. The arc links the end-of-pairing node at the pairing arrival to the midnight node at the end of the next three days.

Arc costs depend on the arc type. Let $i$ and $j$ be two nodes linked by an arc in $A^k$ and let $c_{ij}^k$ be its cost. Rest and postpairing rest arcs may have a positive cost depending on the pilot's preferences, while all other costs are zero. The schedule associated with selecting arcs forming a full path in the network has a cost equal to the sum of the arc costs, which is equal to costs \eqref{eq:costs}. Then, using the dual variables defined above, we set reduced arc costs as follows and where $i$ and $j$ here refer to either the node identifier or the corresponding activity depending on context:
\[
\bar{c}_{ij}^{k} =
\begin{cases}
c_{ij}^{k} - \alpha_{i} & \text{for pairing arcs},\\
c_{ij}^{k} - \beta_{i} & \text{for rest arcs that ensure respecting a preassigned day},\\
c_{ij}^{k} - \gamma^{k} & \text{for roster start arcs, i.e., } i = o(k), \\
c_{ij}^{k} & \text{otherwise}.
\end{cases}
\]

Taking the sum of all reduced costs for the arcs in a schedule gives the same result as Equation \eqref{eq:dual_gen}.

In addition, there are resource windows on each node. Resources are values that are consumed when traversing arcs. Resource windows require that at each node, each resource fall within a given range. These ranges ensure that the schedules generated are feasible. Moreover, the lower bounds on the resource windows are soft: it is possible to extend a schedule even if it violates the lower bound of the arrival node. In this case, the corresponding resource values are set to the lower bound.

We define the resources as follows. The \textit{time} resource is the chronological time, in minutes, elapsed since the beginning of the month. Each assigned arc increases this resource by the amount of time needed to perform the activity encoded by the arc. The \textit{days off} resource is the number of days off left to grant to the pilot to meet the monthly minimum for time off. It is initialized at $T^{off}$, and rest, vacation, postpairing rest and postpairing vacation arcs decrease its value by the number of days off taken. To ensure that pilot $k$ receives enough days off, the resource window on the sink node is $[0, 0]$ days off left. The \textit{flight  time} resource is the flight time cumulated by the pilot so far, in hours. Pairing arcs increase the value of this resource by the number of hours of work included in the flights in the pairing and the briefing and debriefing at the start and end of the pairing. This resource is initialized at 0 at the source node and must respect a resource window of $[0, T^{flight}]$ at each node in the network. The last resource is the \textit{number of consecutive days worked} (i.e., consecutive days where the pilot serves at least one departing flight). It is initialized at 0, and its value increases when selecting pairing arcs by the number of days in the pairing where at least one flight departs, while rest arcs, including postpairing rest arcs, reset it to zero. Other arc types have no effect on this resource, and its range is $[0,T^{work}]$ consecutive working days at every node.

When no new columns are generated by the subproblems, we perform a branching step. However, we note that practical implementations of the CRP do not always solve each linear relaxation to optimality, but rather set a heuristic stopping criterion. In this paper, we require that the objective improve by a percentage of at least $m_{iter}$ over $N_{iter}$ rolling iterations. When this is not the case, a branching step takes place, whether the linear relaxation is solved to optimality using this criterion or not. We use two branching methods: column fixing and intertask fixing.

Column fixing sets a schedule assignment variable to 1 permanently. Intertask fixing imposes or forbids that two activities appear consecutively in a solution. Schedules that do not respect this condition are not included in the next RMP. In each case, the branching we do is heuristic with a depth-first search that continues until an integer solution is found (i.e., there is no backtracking or branching to infeasibility). We start with column fixing of one to three variables whose fractional value is above a certain threshold (\textit{CfixSelectThreshold}). Intertask decisions are scored next based on fractional flow, which is then standardized based on the scores for all such possible decisions. One to three intertask decisions with a score above a minimum threshold (\textit{ItimposeSelectThreshold}) are then selected. 

\subsection{Windowing method} \label{win_adapt}
We now describe the main approach used to leverage ML in this paper: windowing. When applying windowing, the horizon is broken into several overlapping windows. For example, if windows are taken to be 10 days in length with an overlap of three days, and the month chosen has 31 days, we would have that the windows range from
\begin{itemize}
    \item day 1 - day 10,
    \item day 8 - day 17,
    \item day 15 - day 24, and
    \item day 22 - day 31.
\end{itemize}

Once these windows have been defined, windowing can proceed in two ways. If no initial solution has been provided as input to the windowing procedure, we build a solution from scratch. The solution is built one window at a time and in sequence while holding fixed the elements found in previous windows. We call this basic windowing, \textit{win-basic} for short as used in Section \ref{method}. Basic windowing is the approach that has been used so far in the literature, such as in the papers cited in Section \ref{tradition}.

In the case where an initial solution has been provided as input, the whole solution is reoptimized one window at a time. Windowing improves each part of the roster successively and in sequence, while the elements included in the other windows and outside the overlap period remain fixed. We call this windowing with ML, \textit{win-ML} for short as used in Section \ref{method}. In order to implement both forms of windowing (i.e., \textit{win-basic} and \textit{win-ML}), we solve the CRP with the same branch-and-price scheme that would otherwise be standard, with the difference that the subproblem network for each pilot is modified to impose that the content of some windows remain fixed. 

With \textit{win-basic} (no initial solution), no pairing can be assigned prematurely in future windows. The subproblem networks for the pilots are therefore changed to account for this, and the arcs that would allow assigning pairings starting in future windows, unlike rest arcs, are removed. Furthermore, pairings assigned in previous windowing steps must be selected again (unless the departure occurs in the overlapping period between the previous and current windows), and arcs that are only in paths allowing to bypass this restriction are excluded as well.

With \textit{win-ML} (an initial solution is provided), pairings assigned both in a previous windowing step (with the exception of overlap as above) or future windows of the initial solution, must be respected. For this reason, in each subproblem network, arcs that are in conflict with the pairings in the fixed windows are removed from the network, forcing the solver to select pairings departing in them. We also note that for both types of windowing, nodes that are not part of any allowed path are removed.

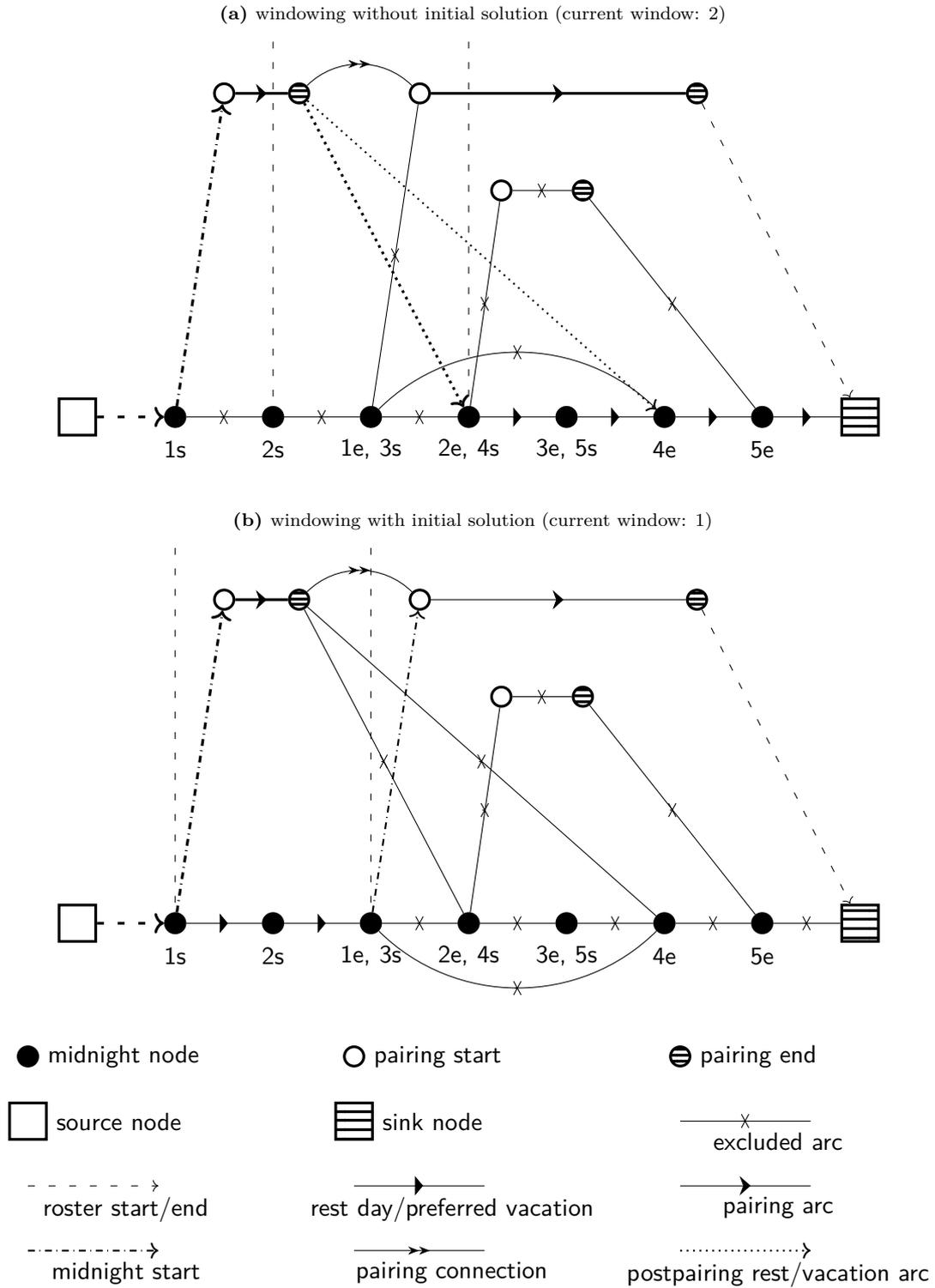
\begin{figure}[ht!]
\centering
\captionsetup{justification=centering, margin=1.75cm}
\begin{minipage}{\linewidth}
\centering

\begin{subfigure}{\textwidth}
\centering
\subcaption{windowing without initial solution (current window: 2)}

\begin{tikzpicture}
[start/.style={circle,draw=black, very thick, inner sep = 3pt},
finish/.style={circle,draw=black, very thick, inner sep = 3pt, pattern={Lines[angle=135,distance=2.5pt]}},
midnight/.style={circle,draw=black, very thick, fill=black, inner sep = 3pt},
source/.style={rectangle,draw=black, very thick, inner sep = 8pt},
sink/.style={rectangle, draw=black, very thick, inner sep = 8pt, pattern={Lines[angle=135,distance=4pt]}}]

\node at (-6,0) [source] (source0) {};
\node at (-3.75, 5) [start] (start01) {};
\node at (-4.5,0) [midnight] (midnight01) {};
\node at (-2.6, 5) [finish] (finish01) {};
\node at (-3,0) [midnight] (midnight02) {};
\node at (-1.5,0) [midnight] (midnight03) {};
\node at (0,0) [midnight] (midnight04) {};
\node at (-0.75, 5) [start] (start02) {};
\node at (1.5,0) [midnight] (midnight05) {};
\node at (3.5, 5) [finish] (finish02) {};
\node at (3,0) [midnight] (midnight06) {};
\node at (0.5,3.5) [start] (start03) {};
\node at (4.5,0) [midnight] (midnight07) {};
\node at (1.75,3.5) [finish] (finish03) {};
\node at (6,0) [sink] (sink0) {};

\node at (-4.5, -0.5) {1s} ;
\node at (-3, -0.5) {2s} ;
\node at (-1.5, -0.5) {1e, 3s} ;
\node at (0, -0.5) {2e, 4s} ;
\node at (1.5, -0.5) {3e, 5s} ;
\node at (3, -0.5) {4e} ;
\node at (4.5, -0.5) {5e} ;

\draw [->, loosely dashed, very thick] (source0) -- (midnight01);
\draw (midnight01) -- node {\crossarrow} (midnight02);
\draw (midnight02) -- node {\crossarrow} (midnight03);
\draw (midnight03) -- node {\crossarrow} (midnight04);
\draw (midnight04) -- node {\midarrow}  (midnight05);
\draw (midnight05) -- node {\midarrow}  (midnight06);
\draw (midnight06) -- node {\midarrow} (midnight07);
\draw (midnight07) -- node {\midarrow} (sink0);

\draw[->, dash dot, very thick] (midnight01) -- (start01);
\draw[very thick] (start01) -- node {\whitearrow} (finish01);
\draw[->, dotted, very thick] (finish01) -- (midnight04);
\draw[->, dotted, thick] (finish01) -- (midnight06);

\draw (midnight03) -- node {\crossarrow}  (start02);
\draw[very thick] (start02) -- node {\whitearrow} (finish02);
\draw [->, loosely dashed] (finish02) -- (sink0);

\draw(midnight04) -- node {\crossarrow}  (start03);
\draw (start03) -- node {\crossarrow} (finish03);
\draw (finish03) -- node {\crossarrow} (midnight07);

\draw (finish01) to[out=45, in=135] node {\forwardarrow} (start02);
\draw (midnight03) to[out=45, in=135] node {\crossarrow} (midnight06);

\draw [loosely dashed] (-3,0) -- (-3,6);
\draw [loosely dashed] (0,0) -- (0,6);

\end{tikzpicture}
\label{fig:win-noinit}
\end{subfigure}

\bigskip

\begin{subfigure}{\textwidth}
\centering
\subcaption{windowing with initial solution (current window: 1)}

\begin{tikzpicture}
[start/.style={circle,draw=black, very thick, inner sep = 3pt},
finish/.style={circle,draw=black, very thick, inner sep = 3pt, pattern={Lines[angle=135,distance=2.5pt]}},
midnight/.style={circle,draw=black, very thick, fill=black, inner sep = 3pt},
source/.style={rectangle,draw=black, very thick, inner sep = 8pt},
sink/.style={rectangle, draw=black, very thick, inner sep = 8pt, pattern={Lines[angle=135,distance=4pt]}}]

\node at (-6,0) [source] (source0) {};
\node at (-3.75, 5) [start] (start01) {};
\node at (-4.5,0) [midnight] (midnight01) {};
\node at (-2.6, 5) [finish] (finish01) {};
\node at (-3,0) [midnight] (midnight02) {};
\node at (-1.5,0) [midnight] (midnight03) {};
\node at (0,0) [midnight] (midnight04) {};
\node at (-0.75, 5) [start] (start02) {};
\node at (1.5,0) [midnight] (midnight05) {};
\node at (3.5, 5) [finish] (finish02) {};
\node at (3,0) [midnight] (midnight06) {};
\node at (0.5,3.5) [start] (start03) {};
\node at (4.5,0) [midnight] (midnight07) {};
\node at (1.75,3.5) [finish] (finish03) {};
\node at (6,0) [sink] (sink0) {};

\node at (-4.5, -0.5) {1s} ;
\node at (-3, -0.5) {2s} ;
\node at (-1.5, -0.5) {1e, 3s} ;
\node at (0, -0.5) {2e, 4s} ;
\node at (1.5, -0.5) {3e, 5s} ;
\node at (3, -0.5) {4e} ;
\node at (4.5, -0.5) {5e} ;

\draw [->, loosely dashed, very thick] (source0) -- (midnight01);
\draw (midnight01) -- node {\midarrow} (midnight02);
\draw (midnight02) -- node {\midarrow} (midnight03);
\draw (midnight03) -- node {\crossarrow}  (midnight04);
\draw (midnight04) -- node {\crossarrow}  (midnight05);
\draw (midnight05) -- node {\crossarrow}  (midnight06);
\draw (midnight06) -- node {\crossarrow}  (midnight07);
\draw (midnight07) -- node {\crossarrow}  (sink0);

\draw[->, dash dot, very thick] (midnight01) -- (start01);
\draw[very thick] (start01) -- node {\whitearrow} (finish01);
\draw(finish01) -- node {\crossarrow}  (midnight04);
\draw (finish01) -- node {\crossarrow} (midnight06);

\draw[->, dash dot, thick] (midnight03) -- (start02);
\draw (start02) -- node {\whitearrow} (finish02);
\draw [->, loosely dashed] (finish02) -- (sink0);

\draw (midnight04) -- node {\crossarrow}  (start03);
\draw (start03) -- node {\crossarrow} (finish03);
\draw (finish03) -- node {\crossarrow}  (midnight07);

\draw (finish01) to[out=45, in=135] node {\forwardarrow} (start02);
\draw (midnight03) to[out=-45, in=-135] node {\crossarrow}  (midnight06);

\draw [loosely dashed] (-4.5,0) -- (-4.5,6);
\draw [loosely dashed] (-1.5,0) -- (-1.5,6);

\end{tikzpicture}
\label{fig:win-init}
\end{subfigure}

\bigskip

\begin{tikzpicture}
[start/.style={circle,draw=black, very thick, inner sep = 3pt},
finish/.style={circle,draw=black, very thick, inner sep = 3pt, pattern={Lines[angle=135,distance=2.5pt]}},
midnight/.style={circle,draw=black, very thick, fill=black, inner sep = 3pt},
source/.style={rectangle,draw=black, very thick, inner sep = 8pt},
sink/.style={rectangle, draw=black, very thick, inner sep = 8pt, pattern={Lines[angle=135,distance=4pt]}}]

\node at (-6, -2) [midnight, label=right:midnight node] {};
\node at (-1, -2) [start, label=right:pairing start] {};
\node at (4, -2) [finish, label=right:pairing end] {};
\node at (-6, -3) [source, label=right:source node] {};
\node at (-1, -3) [sink, label=right:sink node] {};
\draw (4, -3) -- node {\crossarrow} (6, -3) node [near end, below=1pt] {excluded arc};


\draw [->, loosely dashed] (-6, -4) -- (-4, -4) node [near end, below=1pt] {roster start/end};
\draw (-1, -4) -- node {\midarrow} (1, -4) node [near end, below=1pt] {rest day/preferred vacation};
\draw (4, -4) -- node {\whitearrow} (6, -4) node [near end, below=1pt] {pairing arc};
\draw [->, dash dot, thick] (-6, -5) -- (-4, -5) node [near end, below=1pt] {midnight start};
\draw (-1,-5) -- node {\forwardarrow} (1,-5) node [near end, below=1pt] {pairing connection};
\draw[->, dotted, thick] (4,-5) -- (6,-5) node [near end, below=1pt] {postpairing rest/vacation arc};


\end{tikzpicture}

\caption{Pilot subproblem networks with frozen subpaths}
\label{fig:frozen}
\end{minipage}
\end{figure}

Figure \ref{fig:frozen} illustrates the windowing method for windows with a length of 2 days and an overlap of 1 day. It represents two modified subproblem networks based on Figure \ref{fig:roster}, one for each form of windowing. A window start or end is represented by its number followed by the letter \textbf{s} or \textbf{e}, respectively (e.g., \textbf{2s} and \textbf{2e} for the second window). In this example, we focus on the second window for windowing without an initial solution (Figure \ref{fig:win-noinit}), and on the first window for windowing with an initial solution (Figure \ref{fig:win-init}).

In the case without an initial solution (Figure \ref{fig:win-noinit}), we consider the first pairing of the month. The network excludes arcs that conflict with this pairing, as that pairing was assigned at the previous windowing step. We allow to choose whether to assign the pairing starting on the third day, as this falls within the second window (i.e., the one under treatment). The network excludes arcs allowing the pairing starting on the fourth day, as later pairings can only be assigned when the window of their departure is reached.

In the case with an initial solution (Figure \ref{fig:win-init}), we consider that the arcs included in the network must allow to choose whether the pairing starting on the first day is assigned, even though it was part of the ML solution, as the first window is currently treated. In the case where we choose to assign the pairing again, the networks will need to take this into account when treating future windows by imposing its selection, just as in Figure \ref{fig:win-noinit}; in the opposite case, it will be forbidden instead. We impose the pairing starting on the third day, as it is part of the initial solution and does not start in the window treated. The pairing starting on the fourth day is not allowed, as it conflicts with the fixed pairing starting on the third day. We exclude arcs that can only be part of paths that do not select the pairing starting on the third day. Overall, the networks in Figures \ref{fig:win-noinit} and \ref{fig:win-init} are in accordance with the description of the network modifications given above: however, they exclude the pairing starting on the fourth day for different reasons. Using modified networks for each window successively leads to a full solution.

A general feature of windowing is that once all appropriate arcs have been removed in a network, its size is significantly reduced. For example, for the instance with the largest networks in our data, the number of remaining arcs for a window length of 10 days with an overlap of 3 days varies between 66462 and 86661 depending on the window and the presence of an initial solution, by comparison with 263377 for an unchanged network. Since the computational complexity of the subproblems increases approximately in proportion with the square of the number of arcs, these smaller networks lead to an accelerated resolution while applying the already heuristic branch-and-price algorithm. Furthermore, the maximum depth reached in the branching trees is also reduced (from 322 levels to at most 63 for the example given).

We close this section with a word on the challenges presented by windowing. As mentioned in Section~\ref{revue}, windowing has been used for a range of applications. However, it has been assumed until now that windowing might not suit the CRP as, by contrast with other problems such as the CPP where windowing approaches have been used successfully, the CRP includes many global constraints that extend for the entire time horizon of one month. In other words, these constraints overlap with all windows, rather than just one or two as when dealing with pairings. It is therefore necessary to implement windowing to see how well a state-of-the-art CRP solver can handle this additional complexity. Doing so also allows testing the possible benefits of an interaction between windowing and ML for the CRP: ideally, a good ML-generated roster helps in managing the added complexity of this problem better than basic windowing. We could imagine that fixing good pairings in later windows prevents assigning conflicting pairings in the current window. Another possibility is that without an initial solution, days off are not taken into account while solving the first window, as only days off are assigned in the other ones. This can lead to sacrificing pilot preferences later on to ensure that the corresponding resource constraints are met in the end.

\section{Experimental considerations}\label{data_method}
In this section, we describe the specific data sets we use for our experiments, as well as the conditions in which these experiments are carried out. We start by describing the data, and continue with a description of some solving parameters and performance indicators.

\subsection{Data} \label{Data}

We now present the data used for the experiments conducted in this paper, which seek to determine if the interaction of ML and windowing has the potential to improve the resolution of the CRP. The data set we used is directly taken from the material where it was originally published \cite{kasirzadeh2017}. The original data set is made of seven instances (labeled I1 to I7), each with a different fleet size. Each fleet's instance is inspired by real data provided by a major airline carrier. We define an instance as a set of flights, a set of pilots, and a pairing solution for the CPP given the same set of flights. Furthermore, each pilot is attached to exactly one base in $D$, with $D$ of size 3. The CPP solutions are for each instance obtained using a state-of-the art CPP solver: see Quesnel, Desaulniers and Soumis \cite{quesnel2017}.

In this work, we focus on the two largest instances by number of pairings (I5 and I7). We provide the detailed information concerning the structure of each instance in Table~\ref{table:data}. The first five columns of the table give the instance and number of pilots, airports, flights, and pairings in the instance, respectively. The last two columns represent the percentage of the pairings with length $\geq 4$ days (i.e., long pairings) for the base with the most long pairings, and the total flight time required for all the pairings in the roster, respectively.

\begin{table}[ht]
\centering
\captionsetup{justification=centering}
\begin{tabular}{rrrrrrr} \toprule
Instance & Pilots & Airports & Flights & Pairings & Long pairings (\%) & Flight time (min) \\ \midrule
I5 & 239 & 34 & 5743 & 1222 & 64.5 & 1037333 \\
I7 & 322 & 54 & 7765 & 1473 & 57.0 & 1289546 \\ \bottomrule
\end{tabular}

\caption{Description of the instances available}
\label{table:data}
\end{table}

Table \ref{table:data} shows that the percentage of long pairings is high for these two instances, which is a factor that increases the complexity of the resolution. Indeed, shorter pairings allow for a greater variety of pairing arrangements and create fewer potential conflicts with pairing patterns that emulate an optimal roster.

A key parameter of the CRP is the average flight time per pilot. Racette \textit{et al.} \cite{racette2024} show that this value correlates with the CRP's difficulty. In other words, instances where average flight time per pilot is higher are more constraining and harder to solve. To test our method under various difficulty conditions, we create new instances by changing the number of pilots. For both I5 and I7, the number of pilots available is chosen so that the average flight time per pilot required for the month varies between 50 and 75 hours (to the nearest hour), by increments of 5 hours.

Moreover, preferences, vacations and preassigned days off are not a priori part of the an instance, and must be created. Several scenarios are created for each new instance. A scenario consists of a new instance plus a set of preferences, i.e., preferred vacations and preferred days off. We generate 30 scenarios for each new instance, of which 25 are used for training and the remaining 5 for testing of the ML model used.

\subsection{Test parameters and indicators} \label{method}
In our tests, we consider five different solution methods. The first is a standard branch-and-price method (without windowing). It was implemented using GENCOL 4.5, i.e., commercial  software that is regularly used in industry to solve complex crew scheduling tasks, including the CRP. We call this method \textit{alg-basic}, and we refer to it as our benchmark for all other methods in terms of speed and solution quality. 

The next method is the \textit{seqAsg} procedure by itself. In this case, the roster produced by sequential assignment is kept as a final solution. The remaining three methods are windowing starting with an ML solution (\textit{win-ML}), basic windowing (\textit{win-basic}), and a fast GENCOL heuristic with modified stopping criteria (\textit{alg-fast}). \textit{win-ML} uses the initial solution provided by \textit{seqAsg} and serves to determine if ML has the potential to enhance basic windowing. \textit{win-basic}, on the other hand, serves a benchmark for performance with the use of windowing.

The fast heuristic (\textit{alg-fast}) modifies \textit{alg-basic} by changing the stopping criteria $N_{iter}$ and $m_{iter}$, which set the number of rolling iterations where a minimum improvement of the objective is required. Resource dominance is also changed so dominance goes down from three resources to one resource only. We did not notice any significant difference in performance based on the choice of the resource kept, although the number of resources considered for dominance does have an impact on computational speed.

We now give the details involved in getting ML-generated solutions. ML training is done using a neural network approximator with two hidden layers. These layers have 4 neurons each. With 13 neurons on the input layer and 1 only as output, we get a total of of 77 ML weights when including bias terms as well. This is smaller than is typically the case for a neural network, as the CMA-ES algorithm generally works better for up to 100 learning parameters \cite{back2013}. The activation function selected is the rectified linear unit (ReLU), where $f:\mathbb{R} \rightarrow \mathbb{R}$, $f(x) = \max(0,x)$. Early stopping is used after 120 iterations to prevent overfitting of the data. Furthermore, the \textit{seqAsg} procedure is implemented in Python version 3.7.1. The assignment problems are implemented using the \textit{pulp} library and solved using the CBC solver.

Table \ref{table:fixed_consts} shows the value for different problem parameters, including the upper bounds for each resource, the penalty weights in objective \eqref{eq:crp:obj}, and heuristic branching parameters. The last two parameters in the table refer to running GENCOL using a fast heuristic with stricter stopping criteria before branching. The required improvement and maximum number of iterations to reach it have been changed by factors of 20 and 125, respectively when compared to alg-basic. We also limit resource dominance to one resource instead of three.

Regarding windowing, we set the window length to 10 days, and the overlap to 3 days. These values were taken as a compromise. Indeed, shorter windows tend to lead to objectives of a worse quality, whereas larger windows take longer to solve.

Let us now define the key performance indicators that help measuring the performance of each method. We define $S^{*}$ and $t^{*}$ as the mean CRP objective and solving time when using a given method $*$. These indicators taken with the \textit{alg-basic} method serve as a baseline for other solution approaches and the five test scenarios used. We also note that the values above include penalties for any unassigned pairings, and that $t^{*}$ includes the time needed by \textit{seqAsg} to generate ML solutions, if applicable.

We also define:

\begin{itemize}
    \item $L^{*} = 100*\left(1-\frac{S^{*}}{S^{alg-basic}}\right)$ is the loss of satisfaction for solutions obtained with a given method $*$ when compared to \textit{alg-basic}.
    \item $p = 100*\frac{t^{*}}{t^{alg-basic}}$ is the proportion of the solving time taken by a given method $*$ when compared to \textit{alg-basic}.
\end{itemize}

These performance indicators allow to compare and constrast the use of windowing with both a standard branch-and-price algorithm and two accelerated resolution methods (the fast GENCOL heuristic and stand-alone sequential assignment). The differentiation between \textit{win-ML} and \textit{win-basic} helps clarify if any effect observed when using windowing can be improved by ML, i.e., if the nature of the input given has an impact.

\begin{table}[H]
\centering
\captionsetup{justification=centering}
\scalebox{0.75}{
\begin{tabular}{rrp{7cm}} \toprule
\textbf{Constant} & \textbf{Value} & \textbf{Description}\\ \midrule
\textbf{Resource bounds} & & \\
$T^{work}$ & 6 & Max. number of consecutive duties\\
$T^{off}$ & 10 & Min. number of days off \\
$T^{min}$ & 12 & Min. postpairing rest time (hours) \\
$T^{flight}$ & 85 & Max. flight time (hours) per pilot \\
\\
\textbf{Penalties} & & \\
$C^{F}$ & 100 & Penalty for unassigned flight \\
$C^{D}$ & 1,000,000 & Penalty for preassigned days off \\
\\
\textbf{Resolution} & & \\
$N_{iter}$ & 500 & Number of iter. for min. improvement\\
$m_{iter}$ (\%) & 0.05 & Min. improvement over $N_{iter}$ \\
CfixSelectThreshold & 0.70 & Column fixing threshold \\
ItimposeSelectThreshold & 0.70 & Intertask fixing threshold \\
\\
\textbf{Fast heuristic} & & \\
$N_{iter}$ & 4 & Number of iter. for min. improvement\\
$m_{iter}$ (\%) & 1.00 & Min. improvement over $N_{iter}$ \\

\bottomrule
\end{tabular}
}

\caption{Solving constants by category}
\label{table:fixed_consts}
\end{table}

\section{Computational results} \label{results}
In this section, we present results for the various solution methods outlined in this paper. In Section \ref{win_alg}, we first compare the quality and time needed for solutions obtained with windowing to those produced by a standard branch-and-price algorithm. In Section \ref{win_ml}, we compare the windowing method with (\textit{win-ML}) and without (\textit{win-basic}) an initial solution. Finally, in Section \ref{win_len}, we study the impact of changing the length of the windows on the algorithm's performance.

\subsection{Impact of windowing} \label{win_alg}
In this section, we compare the speed and objective value for solutions obtained with windowing alone (\textit{win-basic}) to those obtained by the standard branch-and-price algorithm (\textit{alg-basic}) and the accelerated GENCOL heuristic \textit{alg-fast}.

We report these results in Table \ref{table:win_alg}. Each row represents a different instance with a given level of average total flight time per pilot, rounded to the nearest hour. For each instance, the results are categorized based on the method that was used to obtain them, i.e., \textit{alg-basic}, \textit{alg-fast}, or \textit{win-basic}. For each method, we report the average loss ($L$) as well as the average computing time, in seconds ($t$). For \textit{alg-fast} and \textit{win-basic}, we also report the computing time relative to \textit{alg-basic} ($p$).

\begin{table}[ht!]
\centering
\captionsetup{justification=centering}
\scalebox{0.90}{

\begin{tabular}{\$r^r^r^r^r^r^r^r^r^r^r^r} \toprule
  \multirow{2}{*}{\textbf{Data}} & \multirow{2}{*}{\textbf{Pilots}} & \multicolumn{2}{c}{\textbf{alg-basic}} & \multicolumn{4}{c}{\textbf{alg-fast}} & \multicolumn{4}{c}{\textbf{win-basic}} \\  \cmidrule(lr){3-4} \cmidrule(lr){5-8} \cmidrule(lr){9-12}
 & & $S$ & $t$ (s) & $S$ & $t$ (s) & $L$ & $p$ (\%) & $S$ & $t$ (s) & $L$ & $p$ (\%) \\ \midrule
 I5-50 & 346 & 253757 & 391 & 250136 & 223 & 1.43 & 57.0 & 252777 & 27.7 & 0.39 & 7.1 \\
 I5-55 & 315 & 231875 & 405 & 229165 & 225 & 1.17 & 55.6 & 230972 & 30.7 & 0.39 & 7.6 \\
 I5-60 & 288 & 210863 & 424 & 208176 & 223 & 1.27 & 52.6 & 209643 & 31.1 & 0.58 & 7.3 \\
 I5-65 & 266 & 195573 & 549 & 192070 & 246 & 1.79 & 44.8 & 193791 & 35.4 & 0.91 & 6.4 \\
 I5-70 & 247 & 181200 & 691 & 176674 & 244 & 2.50 & 35.3 & 178209 & 39.3 & 1.65 & 5.7 \\
 I5-75 & 230 & 167501 & 769 & 163200 & 243 & 2.27 & 31.6 & 163162 & 44.0 & 2.59 & 5.7 \\
 \\
 I7-50 & 429 & 316174 & 633 & 311408 & 290 & 1.51 & 45.8 & 314419 & 35.1 & 0.51 & 5.5 \\
 I7-55 & 391 & 288098 & 680 & 284343 & 278 & 1.30 & 40.9 & 286515 & 44.6 & 0.55 & 6.6 \\
 I7-60 & 358 & 264681 & 790 & 259883 & 263 & 1.81 & 33.3 & 261668 & 50.8 & 1.14 & 6.4 \\
 I7-65 & 331 & 245747 & 898 & 240191 & 229 & 2.26 & 25.5 & 241898 & 52.6 & 1.57 & 5.9 \\
 I7-70 & 307 & 225507 & 1228 & 218704 & 227 & 3.02 & 18.5 & 220418 & 64.5 & 2.26 & 5.3 \\
 I7-75 & 287 & 207255 & 1477 & 200403 & 220 & 3.31 & 14.9 & 199769 & 69.7 & 3.61 & 4.7 \\
 \\
 \textbf{Average} & \textbf{316} & \textbf{232353} & \textbf{745} & \textbf{227863} & \textbf{243} & \textbf{1.97} & \textbf{38.0} & \textbf{229437} & \textbf{43.8} & \textbf{1.35} & \textbf{6.2} \\
\bottomrule

\end{tabular}}

\caption{Quality and solving time of solutions obtained with \textit{win-basic} compared to \textit{alg-basic} and \textit{alg-fast}} 
\label{table:win_alg}
\end{table}

From Table \ref{table:win_alg}, we observe that \textit{win-basic} exhibits computing times sped up on average by a factor greater than 10 with respect to \textit{alg-basic}. In fact, in many cases, we have that $p^{win-basic}$ is closer to 5\% (e.g., for I7 and $W=75$). Furthermore, the values of $p^{win-basic}$ show that the solving time grows more slowly when using windowing than for standard resolution when the CRP problem instances become progressively more difficult. For example, in all cases, we have that the the reported value for $p^{win-basic}$ is lower when $W \in \{70,75\}$ than when $W \le 65$.

Method \textit{win-basic} consistently performs faster than \textit{alg-fast}, with a better average objective value. This is significant, as since this heuristic was configured to find a solution as quickly as possible, we have evidence that naive heuristics cannot outperform windowing in this context.

\subsection{Impact of ML solutions} \label{win_ml}
In this section, we test the impact of providing the windowing procedure with ML input (\textit{win-ML}). We compare \textit{win-ML} with basic windowing (\textit{win-basic}) and the use of ML solutions alone (\textit{seqAsg}).

We show the results for these comparisons in Table \ref{table:ml_tests}.

\begin{table}[ht!]
\centering
\captionsetup{justification=centering}
\scalebox{0.85}{

\begin{tabular}{\$r^r^r^r^r^r^r^r^r^r^r^r^r^r} \toprule
  \multirow{2}{*}{\textbf{Data}} & \multirow{2}{*}{\textbf{Pilots}} & \multicolumn{2}{c}{\textbf{seqAsg}}  & \multicolumn{4}{c}{\textbf{win-basic}} & \multicolumn{4}{c}{\textbf{win-ML}} \\  \cmidrule(lr){3-4} \cmidrule(lr){5-8} \cmidrule(lr){9-12}
 & & $S$ & $L$ & $S$ & $t$ (s) & $L$ & $p$ (\%) & $S$ & $t$ (s) & $L$ & $p$ (\%) \\ \midrule
 I5-50 & 346 & 236462 & 6.82 & 252777 & 27.7 & 0.39 & 7.1 & 252823 & 40.2 & 0.37 & 10.3 \\
 I5-55 & 315 & 213577 & 7.89 & 230972 & 30.7 & 0.39 & 7.6 & 231105 & 42.6 & 0.33 & 10.5 \\
 I5-60 & 288 & 193213 & 8.37 & 209643 & 31.1 & 0.58 & 7.3 & 210014 & 44.8 & 0.40 & 10.6 \\
 I5-65 & 266 & 185979 & 4.91 & 193791 & 35.4 & 0.91 & 6.4 & 194658 & 48.5 & 0.47 & 8.8 \\
 I5-70 & 247 & 161701 & 10.76 & 178209 & 39.3 & 1.65 & 5.7 & 178740 & 50.7 & 1.35 & 7.3 \\
 I5-75 & 230 & 153856 & 8.15 & 163162 & 44.0 & 2.59 & 5.7 & 164796 & 53.7 & 1.61 & 7.0 \\
 \\
 I7-50 & 429 & 307825 & 2.64 & 314419 & 35.1 & 0.51 & 5.5 & 315346 & 53.9 & 0.26 & 8.5 \\
 I7-55 & 391 & 275503 & 4.37 & 286515 & 44.6 & 0.55 & 6.6 & 287093 & 62.6 & 0.35 & 9.2 \\
 I7-60 & 358 & 252426 & 4.63 & 261668 & 50.8 & 1.14 & 6.4 & 263156 & 64.0 & 0.58 & 8.1 \\
 I7-65 & 331 & 228589 & 6.98 & 241898 & 52.6 & 1.57 & 5.9 & 243622 & 68.5 & 0.86 & 7.6 \\
 I7-70 & 307 & 205510 & 8.87 & 220418 & 64.5 & 2.26 & 5.3 & 221864 & 71.3 & 1.62 & 5.8 \\
 I7-75 & 287 & 182271 & 12.05 & 199769 & 69.7 & 3.61 & 4.7 & 201816 & 85.4 & 2.62 & 5.8 \\
 \\
 \textbf{Average} & \textbf{316} & \textbf{216409} & \textbf{7.20} & \textbf{229437} & \textbf{43.8} & \textbf{1.35} & \textbf{6.2} & \textbf{230419} & \textbf{57.2} & \textbf{0.90} & \textbf{8.3} \\
\bottomrule

\end{tabular}}

\caption{Quality and solving time of solutions obtained with \textit{win-ML} compared to \textit{win-basic} and \textit{alg-ML}}
\label{table:ml_tests}
\end{table}

From the results presented in Table \ref{table:ml_tests}, a number of facts stand out. First, both windowing methods provide much improved objective values when compared to \textit{seqAsg}. In particular, we see that the average values of $L^{seqAsg}$ and $L^{win-ML}$ are quite different, with a drop from 7.20\% to 0.90\%, or a recovery of 87.5\% of the objective loss with respect to $S^{win-basic}$. Next, of primary importance, we can observe that for all instances, $L^{win-ML}$ is lower than $L^{win-basic}$, making \textit{win-ML} better performing than \textit{win-basic}. In other words, we provide a new way of generating accelerated solutions using both ML and windowing that beats windowing alone.

Furthermore, the enhanced quality of the \textit{win-ML} solutions is more marked as the total flight time required from pilots increases. For example, for instance I5 and $W=50$, the value of $L^{win-ML}$ is 0.37\% for the improved ML solution while $L^{win-basic}$ is 0.39 \%.  For $W=75$, these values are 1.61\% and 2.59\%, respectively. Given that values of $W$ near 75 hours flown per month are similar to industry standards and typical requirements placed on pilots, this means our method is particularly efficient in situations presenting a realistic challenge.

The average results for \textit{win-ML} showed in Table \ref{table:ml_tests} also confirm the findings discussed so far. We see that $L^{win-basic}$ is higher than $L^{win-ML}$: the loss of satisfaction with respect to $S^{alg-basic}$ drops from 1.35\% to 0.90\% (a recovery of 33.3\% of the loss observed with basic windowing). Given that the linear relaxation found by \textit{alg-basic} is in every case fully solved due to stringent stopping criteria ($N_{iter}=500$) and that the integrality gap of \textit{win-ML} solutions is under 0.10\%, \textit{win-ML} solutions are on average less than 1\% away from optimal solutions: if this threshold is judged acceptable, then the solutions can be used as final products. These findings are especially interesting in light of the fact that, for each window, the CRP is solved from scratch. This means that the information contained in a window cannot be used during the CRP run aiming to enhance it. This leads to an important question related to the factors that lead to ML-generated solutions performing better than solutions obtained through other means. In particular, good assignments within a window are found independently again after solving the problem for that window: they do not simply remain in place, while the solver leverages them. One hypothesis to explain the observed gains in performance is that since all but the current window remain frozen, the resolution of the CRP for the current window cannot cause the assignment of activities that are incompatible with good activities in the fixed parts of the roster.

Finally, we also see that while \textit{win-ML} takes slightly longer to run than \textit{win-basic}, mostly due to the requirement of generating an ML solution with \textit{seqAsg}, the average running time is still less than 10\% of what is needed in GENCOL with \textit{alg-basic}, as shown by the average value of $p^{win-ML}$ for this method.

\subsection{Impact of window length} \label{win_len}
An important point of consideration is the influence of window length. Short windows generally lead to a worse objective, and given that 10-day windows are already fast, there is little interest in gaining a few seconds of computational time at the cost of a worsened value for the objective function. However, longer windows may improve the objective, and it is for this reason interesting to see if they can do so with a reasonable acceleration with regard to \textit{alg-basic.}

\begin{table}[ht!]
\centering
\captionsetup{justification=centering}

\scalebox{1.00}{

\begin{tabular}{\$r^r^r^r^r^r^r^r^r^r^r} \toprule
  \multirow{2}{*}{\textbf{Data}} & \multicolumn{2}{c}{\textbf{alg-basic}} & \multicolumn{4}{c}{\textbf{win-ML (10 days)}} & \multicolumn{4}{c}{\textbf{win-ML (15 days)}} \\  \cmidrule(lr){2-3} \cmidrule(lr){4-7} \cmidrule(lr){8-11}
 & $S$ & $t$ (s) & $S$ & $t$ (s) & $L$ & $p$ (\%) & $S$ & $t$ (s) & $L$ & $p$  (\%) \\ \midrule
 I5-65 & 195573 & 549 & 194658 & 48.5 & 0.47 & 8.8 & 195350 & 116 & 0.11 & 21.1 \\
 I5-70 & 181200 & 691 & 178740 & 50.7 & 1.35 & 7.3 & 180404 & 126 & 0.44 & 18.2 \\
 I5-75 & 167501 & 769 & 164796 & 53.7 & 1.61 & 7.0 & 166790 & 130 & 0.43 & 16.9 \\
 \\
 I7-65 & 245747 & 898 & 243622 & 68.5 & 0.86 & 7.6 & 245428 & 188 & 0.13 & 21.0 \\
 I7-70 & 225507 & 1228 & 221864 & 71.3 & 1.62 & 5.8 & 224952 & 197 & 0.25 & 16.0 \\
 I7-75 & 207255 & 1477 & 201816 & 85.4 & 2.62 & 5.8 & 205906 & 230 & 0.65 & 15.5 \\
\\
 \textbf{Average} & \textbf{203797} & \textbf{935} & \textbf{200916} & \textbf{63.0} & \textbf{1.42} & \textbf{7.1} & \textbf{203138} & \textbf{165} & \textbf{0.34} & \textbf{18.1} \\
\bottomrule

\end{tabular}}

\caption{Quality of solutions obtained through windowing with ML: 10 and 15-day windows}
\label{table:other_winds}
\end{table}

In Table \ref{table:other_winds}, we use window lengths of 15 days with an overlap of 7 days, for a total of three windows over the entire horizon. We see how \textit{alg-basic} and \textit{win-ML} compare both for 10-day and 15-day windows. The rows and columns included have the same meaning as for Tables \ref{table:win_alg} and \ref{table:ml_tests}. We only keep here instances where $W \ge 65$ to test our method on the most complicated cases (i.e., those with a value of $L^{win-basic}$ near or above 1\%). In doing so, we test whether using a longer window length is advisable.

Based on the data in Table \ref{table:other_winds}, it seems that longer time windows yield on average a better objective quality when using \textit{win-ML}, as $L^{win-ML}$ goes down from 1.42\% to 0.34\%. The value of $p^{win-ML}$, however, increases on average from 7.1 \% to 18.1\% of the time needed by \textit{alg-basic}. Which window length is more desirable is then a matter of the specific needs of the airline: it may be that the need to re-solve the CRP is pressing enough that faster resolution is preferred, or that satisfaction takes precedence. The window length can then be adjusted accordingly.

\section{Discussion} \label{discussion}
In this paper, we have presented a new method to obtain accelerated, high-quality solutions to the CRP. We have considered two large instances adapted from real-life data obtained from a major commercial carrier. We have considered different variations based on these instances to get a broader range of monthly flying times per pilot representative of the industry.

We used the \textit{seqAsg} procedure presented in Racette \textit{et al.} \cite{racette2024}, a sequential assignment procedure that creates ML-generated solutions for the CRP. Once ML solutions were obtained, we improved them using windowing, a technique that solves the CRP again by considering certains parts of a solution and freezing the rest. To do so, we have created the first implementation of this approach for the CRP. Based on the results presented, we were able to produce high-quality solutions in under 10\% of the time that a state-of-the-art solver, GENCOL version 4.5, would otherwise have taken. The loss of objective quality associated with these rosters when compared to optimality was under 1\% on average for the cases considered; with longer windows, we obtained better objectives, but at an increased time cost.

By doing so, we have provided two original contributions. First, we showed that windowing can be successfully implemented for the CRP. Second, we generated high-quality solutions for the CRP in less than one minute on realistic instances. This is new, as while many heuristics, ML techniques and other algorithms have improved the resolution of this problem, the extent to which we have accelerated the process with the \textit{win-ML} method is unprecedented. This of course came at the cost of a slightly higher loss of pilot satisfaction, but our method compared favorably with three other fast heuristics: ML alone (\textit{seqAsg}), basic windowing without an initial solution (\textit{win-basic}), and accelerated branching in GENCOL, which suggests that in cases where speed is critical, an approach merging ML and windowing outperforms many naive methods in terms of the quality of the rosters. The solutions obtained can also be used as final products when an average gap of under 1\% is acceptable with respect to optimality.

Finally, with the work presented in this paper, we have given additional evidence that simple ML methods can be rich in possibilities to treat difficult scheduling tasks such as the CRP. In fact, it is noteworthy that such results could be obtained with our method: a well-established evolutionary algorithm (CMA-ES) with a small neural network. We anticipate that future research in this area will help clarify the potential of different ML/OR interactions on a variety of tasks.

\section*{Acknowledgements}
The authors wish to acknowledge the help of François Lessard for assistance in the technical implementation of the windowing procedure in GENCOL.

\bibliographystyle{adapted-bmc.bst}
\bibliography{main}

\end{document}